\documentclass[12pt]{article}

\usepackage[english]{babel}
\usepackage{amsmath}
\usepackage{amssymb,amsfonts,amsthm}
\usepackage{tikz,color,pgf,graphicx,url}
\usetikzlibrary{calc}
\usetikzlibrary{decorations.pathreplacing}
\usepackage{enumerate}     % več možnosti za sezname

\usepackage{xcolor} % for colors
\definecolor{mygreen}{RGB}{55, 184, 2}

\usepackage{hyperref}
\hypersetup{
	colorlinks   = true, %Colours links instead of ugly boxes
	urlcolor     = black, %Colour for external hyperlinks
	linkcolor    = black, %Colour of internal links
	citecolor   = black %Colour of citations
}

\usepackage[
top=3cm,
bottom=3cm,
inner=2.5cm,      % margini za dvostransko tiskanje
outer=2.5cm,
footskip=40pt     % pozicija številke strani
]{geometry}

\newtheorem{theorem}{Theorem}[section]
\newtheorem{corollary}[theorem]{Corollary}
\newtheorem{lemma}[theorem]{Lemma}
\newtheorem{proposition}[theorem]{Proposition}

\newtheorem{definition}[theorem]{Definition}
\newtheorem{remark}[theorem]{Remark}
\newtheorem{example}[theorem]{Example}

\DeclareMathOperator{\capt}{capt}
\DeclareMathOperator{\diam}{diam}
\DeclareMathOperator{\ecc}{ecc}

\def\N{\mathbb{N}}

\begin{document}
	
	\title{The capture time in the game of cops and many robbers}
	
	\author{
		Miha Gyergyek $^{1,}$\thanks{Email: \texttt{mg81805@student.uni-lj.si}}
		\and
		Vesna Iršič Chenoweth $^{1,2,}$\thanks{Email: \texttt{vesna.irsic@fmf.uni-lj.si}}
	}
	
	\maketitle
	
	\begin{center}
		$^1$ Faculty of Mathematics and Physics,  University of Ljubljana, Slovenia\\
		\medskip
		
		$^2$ Institute of Mathematics, Physics and Mechanics, Ljubljana, Slovenia\\
		\medskip
	\end{center}
	
	\begin{abstract}
		The game of cops and robber is a pursuit-evasion game played on graphs that has been extensively studied. Traditionally the game is played with one or more cops and only one robber, but in this paper we consider the game played with multiple robbers, with rules allowing for consecutive capture of the robbers. While the cop number remains the same, the capture time can differ drastically compared to the game played with one robber. In this paper we present general upper bounds for the capture time of cop-win graphs in the game with multiple robbers and show that they are the best possible for some families of graphs. Moreover, we explore how the capture time behaves with an increasing number of robbers and prove a surprising relation between this behavior and the 0-visibility cop number.
	\end{abstract}
	
	
	\noindent
	{\bf AMS Subj.\ Class.\ (2020):} 05C57; 05C05; 05C75
	
	%%%%%%%%%%%%%%%%%%%%%%%%%%%%%%%%%
	%%%%%%%%%%%%%%%%%%%%%%%%%%%%%%%%%
	\section{Introduction}
	\label{sec:intro}
	{\em The game of cops and robbers} is a pursuit-evasion game played on undirected, reflexive, and connected graphs. It is played by two types of players, {\em the cops} and {\em the robbers}. All players are allowed to move between adjacent vertices in each round. The goal of the cops is to capture the robbers, while the robbers are trying to evade capture.
	
	Traditionally, the game is played by a single robber. Graphs on which $k$ cops always capture the robber are called {\em $k$-cop-win graphs} and the smallest such $k$ is the {\em cop number} $c(G)$ of the graph $G$. The case of a single cop playing against a single robber was first described independently in \cite{Qui78} and \cite{NoWi83}, both of which characterized 1-cop-win graphs (cop-win graphs for short). In the most common version of the game, first looked at in \cite{AiFr84}, multiple cops try to capture one robber.
	
	The cops also want to win as quickly as possible in $k$-cop-win graphs. For this purpose, {\em capture time} was introduced \cite{BONATO20095588}. It is the index of the earliest round in which the cops win, assuming that the robber tries to maximize the number of rounds it takes to capture him.
	
	In another version of the game, called the {\em 0-visibility cops and robbers}, the cops cannot see the robber unless they are on the same vertex as him, while the robber has perfect information, knowing the cops' strategy in advance. This version of the game was introduced in \cite{Tos86}.
	
	In this paper we study the version of the game with multiple robbers which is defined in Section \ref{sec:def}. Adding (a finite number of) robbers does not affect the outcome of the game, however it does affect capture time if all robbers must be caught before the cops win. Note that the game with multiple robbers on directed graphs has briefly been looked at on digraphs in \cite{HAHN20062492}, but does not focus on capture time. In \cite{BONATO20095588}, a set of rules has been introduced for the version of the game with more than one robber, but the paper focuses on only the game with one. Those rules differ from the ones we use in this paper and we compare the two later. Both versions are natural extensions of the original game, but as in our version the cop number of the graph remains the same, we can focus on comparing the change in capture time as the number of robbers grows.
	
	We present upper bounds for the capture time in the game of cops and many robbers on cop-win graphs. 
	We also show that the general upper bound is attained in the game with one cop and three robbers in a subfamily of maximum capture time graphs introduced in \cite{BONATO20095588}. We generalize this further by adding more and more robbers and pose the problem: on which graphs and for which $k$ does the limit $\lim_{m\rightarrow \infty} \capt_k(G,m)$ exist? We explore this limit for several well-known families of graphs and then prove a surprising result that the graphs on which the limit exists are precisely the graphs with cop number $k$ in the 0-visibility game.
	%%%%%%%%%%%%%%%%%%%%%%%%%%%%%%%%%
	%%%%%%%%%%%%%%%%%%%%%%%%%%%%%%%%%

	%%%%%%%%%%%%%%%%%%%%%%%%%%%%%%%%%%%%%%%
	%%%%%%%%%%%%%%%%%%%%%%%%%%%%%%%%%%%%%%%
	\section{Preliminaries}
	\label{sec:prelim}
	
	In this section we first formally define the game with multiple robbers and the generalized capture time, and then recall notation and known results that are needed in the rest of the paper.
	
	\subsection{Rules of the game and the generalized capture time}
	\label{sec:def}
	We generalize the usual game of cops and robbers to include multiple robbers in a slightly different way to \cite{BONATO20095588}.
	We follow the convention used in \cite{BoNo11} that the cops are female and the robbers are male. The set of cops is denoted as $C$, while the set of robbers is denoted as $R$.
	
	The rules are as follows.
	The cops pick their starting vertices first, then the robbers do the same. We call this {\em round 0}. Afterwards, the cops and the robbers may either move to an adjacent vertex or stay on the same vertex. We call this a {\em move}. A \emph{round} consists of all $k$ cops making a move, followed by all $m$ robbers making a move afterwards. The $i$-th round after round 0 is {\em round $i$}.
	
	We allow for multiple players to occupy the same vertex. Whenever a cop occupies the same vertex as some robbers, the robbers 
	at that vertex are all {\em captured} and do not continue playing.
	
	The cops win if all the robbers are captured in a finite number of moves; the robbers win if at least one robber evades capture forever.
	
	Here, if there are at least two robbers, our rules differ from the ones in \cite{BONATO20095588} in the sense that we allow for the robbers to be captured one by one by the cops, while all robbers have to be on a vertex occupied by a cop in the final round in that paper. In the game with only one robber, both definitions align.
	
	We define {\em the winning strategy for the cops} as a sequence of moves (that may depend on the robbers' moves) after which all the robbers are captured despite what the robbers do.
	Similarly, {\em the winning strategy for the robbers} is a sequence of moves (they may again depend on the cops' moves) that allows at least one robber to evade capture forever. As before, the graphs on which $k$ cops always win are {\em $k$-cop-win graphs against $m$ robbers}. For short, we say that 1-cop-win graphs are cop-win graphs.
	
	We immediately note the following.
	
	\begin{lemma}
		\label{lem:win}
		A graph G is a $k$-cop-win graph against $m$ robbers if and only if it is $k$-cop-win against a single robber.
	\end{lemma}
	
	\begin{proof}
		The forward direction is clear. Conversely, if the cops have a strategy that allows them to capture a single robber, they may use this strategy on every robber one after the other to capture them all.
	\end{proof}
	
	Because of this lemma, we will still call $k$-cop-win graphs against $m$ robbers just $k$-cop-win graphs.
	Note that an analogous lemma does not hold for the version of the game from \cite{BONATO20095588} since in that game, at least $m$ cops are needed to catch $m$ robbers if those robbers stay still on distinct vertices on a graph with at least $m$ vertices.
	
	We can now define the capture time in the game of cops and multiple robbers.
	
	\begin{definition}
		Let $G$ be a $k$-cop-win graph and let $\ell \geq k$. {\em The length of the game} is the index $i$ of the round in which the last robber is captured in the game with $\ell$ cops and $m$ robbers. The {\em capture time} is the length of the shortest game which guarantees that the cops win, assuming that the robbers try to evade capture for as long as possible. We denote it as $\capt_\ell(G,m)$ for an integer $m\geq 1$ and omit the index if $\ell=1$.
	\end{definition}

	\subsection{Notation}
	All graphs in this paper are finite, simple, connected and reflexive. 
	
	The distance between two vertices $u$ and $v$ of a graph $G$ is the length of the shortest path between them. We write $d_G(u,v)$ or more commonly $d(u,v)$ when the graph is clear from the context. We also extend this definition to the distance between a vertex and a set of vertices.
	
	\begin{definition}
		For $u\in V(G)$ and $X\subseteq V(G)$, the distance between them is $d(u,X)$, where $d(u,X)=\min\{d(u,x)\mid x\in X\}$.
	\end{definition}
	
	We denote the \emph{eccentricity} of a vertex $v \in V(G)$ by $\ecc(v) = \max \{d(v, x) \mid x \in V(G)\}$, and the \emph{diameter} of $G$ by $\diam(G) = \max\{\ecc(v) \mid v\in V(G)\}$. 
	
	\begin{definition}
		A homomorphism between the graphs $G$ and $H$ is a mapping $f: V(G) \rightarrow V(H)$ for which it holds that if $uv \in E(G)$, then $f(u)f(v) \in E(H)$. We usually simply denote $f: G \rightarrow H$.
	\end{definition}
	
	\begin{definition}
		Let $G$ be a graph and $H$ an induced subgraph of $G$. Then $H$ is called a {\em retract} of $G$  if there exists a homomorphism $f: G \rightarrow H$, for which it holds that $f(x)=x$ for every $x \in V(H)$. In this case $f$ is called a {\em retraction}. If $H=G-u$ for some $u \in V(G)$ is a retract of $G$, we call $H$ a {\em 1-point retract} and the homomorphism $f$ a {\em 1-point-retraction}.
	\end{definition}
	
	We provide notation for some common graph classes. Let $n,m\in \N$. We denote paths on $n$ vertices with $P_n$, complete graphs on $n$ vertices by $K_n$, cycles on $n$ vertices by $C_n$ and complete bipartite graphs with class sizes $m$ and $n$ by $K_{m,n}$. In particular, a graph $K_{1,n}$ is called a \emph{star}. 
	A \emph{caterpillar} is a tree $T$ such that $T$ without all its leaves is a path. A \emph{wheel} on $n$ vertices is the cycle $C_{n-1}$ with all its vertices adjacent to a universal vertex $u$. We denote it $W_n$ and require that $n \geq 4$.
	
\begin{figure}[ht!]
		\centering
		\begin{tikzpicture}[node/.style ={draw, circle}]
			\node[node] (1) at (0,0) {};
			\node[node] (2) at (-1,-0.58) {};
			\node[node] (3) at (0,1.15) {};
			\node[node] (4) at (1,-0.58) {};
			\node[node] (5) at (-2,-1.15) {};
			\node[node] (6) at (0,2.31) {};
			\node[node] (7) at (2,-1.15) {};
			\draw (1) -- (2);
			\draw (1) -- (3);
			\draw (1) -- (4);
			\draw (2) -- (5);
			\draw (3) -- (6);
			\draw (4) -- (7);
		\end{tikzpicture}
		\caption{The subdivided star $T_3$.}
		\label{fig:sub_star}
	\end{figure}
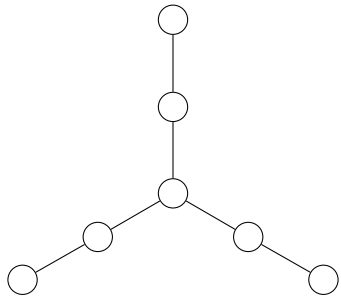

	We define another class of graphs and call it \emph{subdivided stars}. The graph $T_n$ is obtained by taking the star $K_{1,n}$, adding $n$ vertices and connecting each of the added vertices to a distinct leaf of the star. The subdivided star $T_3$ is presented in Figure \ref{fig:sub_star}.
	
	Finally, we recall the following definitions (for example from \cite{BoNo11}).
	A vertex $u\in V(G)$ is a {\em trap} if $N[u]\subseteq N[v]$ for some vertex $v\in V(G)$. We say that $v$ {\em covers} $u$ and write $v\rightarrow u$. These are important in characterizing cop-win graphs.
	
	Let $G$ be a graph, $u \in V(G)$ a trap covered by $v \in V(G)$ and $H$ a 1-point-retract of $G$ with $u$ removed. Then the mapping $f: G \rightarrow H$,
	\begin{equation*}
		f(x) = \begin{cases}
			x; & x\neq u,\\
			v; & x = u
		\end{cases}
	\end{equation*}
	is clearly a retraction which maps a trap to the vertex covering it.
	
	\subsection{The capture time in the game with one robber}
	Let $G$ be a cop-win graph. If time is of the essence for the cops, we may want to know how long it will take them to capture the robber. Specifically, the {\em length} of the game is the round $i$ in which the robber is captured. The robber plays to maximize the length of the game, while the cops try to minimize it. The length of the shortest game under these assumptions is the {\em capture time} $\capt(G)$ of a graph $G$ which was introduced in \cite{BONATO20095588}. We list some results about capture time below.

	\begin{theorem}[{\cite{ClaNow2001}}]
		If $G$ is a cop-win graph on $n \geq 1$ vertices, then $\capt(G) \leq n-1$.
	\end{theorem}
	
	The bound was first improved in \cite{BONATO20095588} for graphs on at least 5 vertices.
	
	\begin{theorem}[{\cite{BONATO20095588}}]
		If $G$ is a cop-win graph on $n \geq 5$ vertices, then $\capt(G)\leq n-3$.
	\end{theorem}
	
	This was further improved by Gavenčiak in \cite{Gav07} for graphs on at least $7$ vertices.
	
	\begin{theorem}[{\cite{Gav07}}]
		\label{thm:n-4}
		If $G$ is a cop-win graph on $n \geq 7$ vertices, then $\capt(G)\leq n-4$.
	\end{theorem}
	The bound in Theorem \ref{thm:n-4} is best possible as is demonstrated by a family of graphs presented in \cite{BONATO20095588} (which is also used later in this paper for a different bound). Additionally, it is proven in \cite{Gav10} that there are at least $2^{n-8}$ non-isomorphic graphs on $n$ vertices with capture time $n-4$ (if $n \geq 8$). However, the structure of these graphs is rather specific and it turns out that a significantly better upper bound can be proven for a large class of graphs that was introduced in \cite{BONATO20095588}. 
	
	\begin{definition}
		A graph $G$ on at least 7 vertices is {\em 2-dismantlable} if it is a cop-win graph and it has two distinct traps $a$ and $b$ and $G-\{a,b\}$ also has two distinct traps or fewer than 7 vertices.     
	\end{definition}
	
	For example, trees and complete graphs are 2-dismantlable.
	
	\begin{theorem}[{\cite{BONATO20095588}}]
		If $G$ is a 2-dismantlable graph on $n$ vertices, then $\capt(G)\leq \left\lfloor \frac{n}{2}\right\rfloor$.
	\end{theorem}
	
	\subsection{0-visibility cops and robbers}
	In the 0-visibility cops and robber, the movement rules of the players are the same as in the original game, but the cops cannot see the robber unless they have captured him. The robber however not only sees the positions of the cops at all times but is also \emph{omniscient}, meaning that the robber knows the complete strategy of the cops. This condition is added to avoid the cops winning ``by chance''. The minimum number of cops needed to capture the robber on a graph $G$ is the \emph{0-visibility cop number} $c_0(G)$. This game was first considered in \cite{Tos86}, is well-understood on trees \cite{dereniowski+2015-0-visibility-2} and with respect to the path-width of the graph \cite{dereniowski+2015-0-vis}. Below we recall some basic results about the game that are needed in the rest of the paper.
	
	\begin{theorem}[{\cite{Tos86}}]
		\begin{enumerate}
			\item For $n\geq 1$, $c_0(K_n)=\left\lceil n/2\right\rceil$.
			\item For $n\geq 1$, $c_0(P_n)=1$.
			\item For $n\geq 3$, $c_0(C_n)=2$.
			\item For $m\geq 1$ and $m\leq n$, $c_0(K_{m,n})=m$.
			\item If $G$ is a tree, $c_0(G)=1$ if and only if $G$ is a caterpillar.
		\end{enumerate}
	\end{theorem}
	
	Additionally, using the bound from \cite{Tos86} that $c_0(G) \geq \omega(G)/2$, we obtain the following.
	
	\begin{proposition}
		\label{prop:0-vis-wheels}
		If $n \geq 3$, then $c_0(W_n) = 2$. 
	\end{proposition}
	
	%%%%%%%%%%%%%%%%%%%%%%%%%%%%%%%%%%%%%%%
	%%%%%%%%%%%%%%%%%%%%%%%%%%%%%%%%%%%%%%%
	
	\section{The capture time of cop-win graphs}
	
	In this section we study the generalized capture time of cop-win graphs. We give several tight upper bounds both in general and for 2-dismantlable graphs. We also determine the capture time for paths, complete graphs, wheels and subdivided stars.
	
	\subsection{Upper bounds}
	We first focus on generalizing Theorem \ref{thm:n-4}.
	
	\begin{definition}
		Let $G$ be a cop-win graph and let the game be played by one cop and one robber. We denote the set of vertices from which the cop requires at most $n-4$ rounds to capture the robber as $Z$.
	\end{definition}
	
	The next lemma provides the most trivial upper bound.
	
	\begin{lemma}
		If $G$ is a cop-win graph on $n \geq 7$ vertices, then $\capt(G,m)\leq (m-1)\diam(G)+m(n-4)$.
	\end{lemma}
	
	\begin{proof}
		By Theorem \ref{thm:n-4}, $\capt(G,1) \leq n-4$. In the worst case for the cop, she has to capture all the robbers one by one.
		It takes $n-4$ rounds from an appropriate vertex to capture the first robber. After that, she may require $\diam(G)$ moves to return to $Z$ before she is able to capture the next robber. In total, she may require $(m-1)\diam(G)+m(n-4)$ rounds to capture all $m$ robbers.
	\end{proof}
	
	Of course, this bound is not very useful for large graphs and we have not found graphs that reach it for any non-trivial number of robbers, so we improve it in the following theorem.
		
	\begin{theorem}
		\label{thm:general_bound}
		If $G$ is a cop-win graph on $n \geq 7$ vertices and
		\begin{equation*}
			d_Z=\max\{d(x,Z)\mid x\in N[u]\text{, $u\in V(G)$ is a trap}\},
		\end{equation*}
		then $\capt(G,m)\leq d_Z(m-1)+m(n-4)$.
	\end{theorem}
	
	\begin{proof}
		As before, the cop starts in $Z$. From there, she captures the first robber in the next $n-4$ rounds. In the worst case, she now requires $d_Z$ moves to return to $Z$, from where she can capture the next robber in $n-4$ moves. In total, it may take her $d_Z(m-1)+m(n-4)$ rounds to capture the last robber.
	\end{proof}
	
	The improvement seems marginal. However, we will show that for a subfamily of maximum capture time graphs constructed in \cite{BONATO20095588}, this bound is tight in a game played by a single cop and at most three robbers. We repeat this construction below.
	
\begin{figure}[ht!]
		\centering
		\begin{tikzpicture}[node/.style args={#1 #2}{draw, circle, label=#1:#2}, scale =0.85, transform shape]
			\node[node=below 6] (6) at (-2.25,0) {};
			\node[node=below 2] (2) at (-0.75,0) {};
			\node[node=below 1] (1) at (0.75,0) {};
			\node[node=below 4] (4) at (2.25,0) {};
			\node[node=below 5] (5) at (0,1) {};
			\node[node=right 3] (3) at (0,-1) {};
			\node[node=right 7] (7) at (0,-2) {};
			\draw (6) -- (2);
			\draw (2) -- (1);
			\draw (1) -- (4);
			\draw (5) -- (6);
			\draw (5) -- (2);
			\draw (5) -- (1);
			\draw (5) -- (4);
			\draw (3) -- (6);
			\draw (3) -- (2);
			\draw (3) -- (1);
			\draw (3) -- (4);
			\draw (7) -- (6);
			\draw (7) -- (4);
			\draw (7) -- (3);
		\end{tikzpicture}
		\caption{The graph $H(7)$.}
		\label{fig:h7}
	\end{figure}
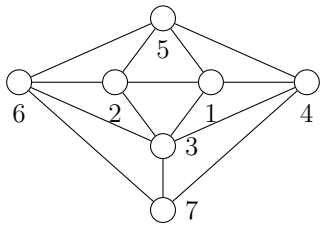

	We begin with the path $P_4$. We add two vertices to the graph, each adjacent to every vertex of the path but not adjacent to each other. We add one last vertex, adjacent to the end-vertices of the path and to one of the two previous vertices, but not to the other vertices. We call this graph $H(7)$; see Figure \ref{fig:h7}. From now on we also use the labels of the vertices as presented in Figure \ref{fig:h7}.
	
\begin{figure}[ht!]
		\centering
		\begin{tikzpicture}[node/.style args={#1 #2}{draw, circle, label=#1:#2}, scale =0.85, transform shape]
			\node[node=below 10] (10) at (-3.75,0) {};
			\node[node=below 6] (6) at (-2.25,0) {};
			\node[node=below 2] (2) at (-0.75,0) {};
			\node[node=below 1] (1) at (0.75,0) {};
			\node[node=below 4] (4) at (2.25,0) {};
			\node[node=below 8] (8) at (3.75,0) {};
			\node[node=below 5] (5) at (0,1) {};
			\node[node=right 9] (9) at (0,2) {};
			\node[node=right 3] (3) at (0,-1) {};
			\node[node=right 7] (7) at (0,-2) {};
			\draw (6) -- (2);
			\draw (2) -- (1);
			\draw (1) -- (4);
			\draw (5) -- (6);
			\draw (5) -- (2);
			\draw (5) -- (1);
			\draw (5) -- (4);
			\draw (3) -- (6);
			\draw (3) -- (2);
			\draw (3) -- (1);
			\draw (3) -- (4);
			\draw (7) -- (6);
			\draw (7) -- (4);
			\draw (7) -- (3);
			\foreach \a/\b in {8/7,8/5,8/4,9/8,9/6,9/5,10/9,10/7,10/6}
			\draw (\a) -- (\b);
		\end{tikzpicture}
		\caption{The graph $H(10)$.}
		\label{fig:h10}
	\end{figure}
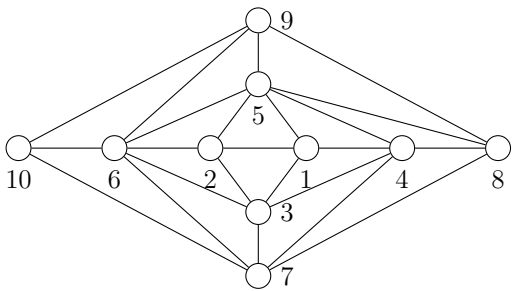

	The graph $H(n)$ for $n \geq 8$ is defined as follows. Take the graph $H(n-1)$, add a vertex $n$ and make it adjacent to vertices $n-1$, $n-3$ and $n-4$. For example, the graph $H(10)$ is presented in Figure \ref{fig:h10}. It was proved in \cite[Theorem 4]{BONATO20095588} that graphs $H(n)$ have maximum capture time, i.e.\ that $\capt(H(n)) = n-4$. For now, we are interested in graphs $H(4\ell+2)$, where $\ell > 1$ is an integer. We prove that graphs $H(4\ell+2)$ attain the bound from Theorem \ref{thm:general_bound} in games with one cop and at most three robbers. We will do this in two steps -- first we will find the set $Z$ for graphs $H(n)$ and then prove that the bound is reached.
	
	\begin{lemma}
		For graphs $H(n)$, $n \geq 7$, the set $Z$ is $\{1, 2\}$.
	\end{lemma}
	
	\begin{proof}
		We denote the cop's vertex immediately after round $i$ by $C^i$ and the robber's vertex by $R^i$.
		
		We first prove that no other vertices are in $Z$. We know from \cite{BONATO20095588} that a robber has a strategy that allows him to survive $n-4$ rounds by starting in vertex $R^0\in [6]$. It is sufficient to show that if the cop does not start the game in either vertex 1 or 2, the robber may survive more than $n-4$ rounds. If the cop starts in $[n] \setminus \{1,2,3,5\}$, the robber may choose such a vertex $R^0\in [6]$ that $d(C^0,R^0)\geq 3$. The robber may then skip his move in the first round and avoid capture for at least $n-4$ more rounds using the strategy from \cite{BONATO20095588}. If the cop does start in either 3 or 5, the robber starts in the other vertex. He then chooses such a vertex $R^1\in [6]$ in round 1 that $d(C^1,R^1)\geq 2$. After that, he is able to keep the game going for at least $n-4$ more rounds, again like in \cite{BONATO20095588}.
		
		We now need to prove that vertices 1 and 2 are in $Z$. The proof for vertex 1 is the same as in \cite{BONATO20095588}. We modify this proof for vertex 2. Assume that the cop starts there. We need to prove that he can capture the robber within the next $n-4$ rounds. Without loss of generality, we may assume that the robber begins in vertex $R^0>2$.
		
		First, suppose that $R^0\not\equiv 0 \pmod 4$. In round 1, the cop moves to such a vertex $C^1$ that $C^1\equiv R^0 \pmod 4$. The cop continues in the following way. Suppose that the robber has not been captured in round $t$ and occupies vertex $R^t$, while the cop occupies $C^t$. Table \ref{tab:tab1} provides the cop's move in round $t+1$. We notice that, since $C^1\equiv R^0 \pmod 4$, it holds that $C^t\equiv R^{t-1} \pmod 4$ for each round $t\geq1$ by induction. Furthermore, $C^t<R^t$ in all rounds except the last round and we notice from the cop's moves in Table \ref{tab:tab1} that once $C^t\equiv R^{t-1} \pmod 4$, the robber cannot increase the difference $R^t-C^t$ in subsequent rounds with any of his moves, so $R^t-C^t$ decreases monotonically. Therefore, the robber will be captured once he is occupying the trap $n$ and the cop is occupying $n-4$. This will happen after $n-5$ rounds at the latest, ending the game in round $n-4$.
		\begin{table}[ht!]
			\centering
			\begin{tabular}{|c|c|}
				\hline
				Robber's vertex after round $t$ & Cop's move in round $t+1$ \\ \hline
				$R^{t-1}-4$         & $C^t+4$                   \\ \hline
				$R^{t-1}-3$         & $C^t+1$                   \\ \hline
				$R^{t-1}-1$         & $C^t+3$                   \\ \hline
				$R^{t-1}$           & $C^t+4$                   \\ \hline
				$R^{t-1}+1$         & $C^t+1$                   \\ \hline
				$R^{t-1}+3$         & $C^t+3$                   \\ \hline
				$R^{t-1}+4$         & $C^t+4$                   \\ \hline
			\end{tabular}
			\caption{Cop's move in relation to the robber's in the first case.}
			\label{tab:tab1}
		\end{table}
		
		The last remaining case to prove is $R^0\equiv 0 \pmod 4$. In this case, the cop's first move will be to $C^1=1$. In round 2, the cop moves in accordance with Table \ref{tab:tab2}.
		\begin{table}[ht!]
			\centering
			\begin{tabular}{|c|c|}
				\hline
				Robber's vertex after round 1 & Cop's move in round 2 \\ \hline
				$R^0-4$         & 4                     \\ \hline
				$R^0-3$         & 5                     \\ \hline
				$R^0-1$         & 3                     \\ \hline
				$R^0$           & 4                     \\ \hline
				$R^0+1$         & 5                     \\ \hline
				$R^0+3$         & 3                     \\ \hline
				$R^0+4$         & 4                     \\ \hline
			\end{tabular}
			\caption{Cop's move in relation to the robber's in round 2 in the second case.}
			\label{tab:tab2}
		\end{table}
		After round 2, the cop proceeds as in the first case. The induction from before now holds for all $t\geq2$. In addition, we notice that $C^2>C^1+1$, so the length of the game doesn't exceed $n-4$. Thus, $2\in Z$.
	\end{proof}
	
	\begin{theorem}
		\label{thm:tight_hn}
		For $n=4\ell+2$ and $\ell\geq2$, the capture time $\capt(H(n),3)$ is $2\ell+3(n-4)$.
	\end{theorem}
	
	\begin{proof}
		We see that since $Z=\{1,2\}$, the distance $d_Z$ is $\ell$.
		
		We will again use the fact from \cite{BONATO20095588} that a robber in vertex $i\in [6]$, who is not adjacent to the cop, can avoid capture for at least $n-4$ more rounds. We denote the robbers $R_1$, $R_2$ and $R_3$ and we find a strategy for one of them to remain uncaptured until at least round $2\ell+3(n-4)$.
		
		The robbers all begin the game in vertex $i\in [6]$ that is not adjacent to the cop's starting vertex. They all play the same strategy that guarantees that none of them is captured before they reach vertex $n$ in round $n-5$ at the earliest. Once they are forced into the trap, the robbers split up between vertices $n$, $n-1$ and $n-3$. The cop may now capture at most one robber in the next round. The remaining robbers now play in the following way. If they are occupying vertex $i$ and the cop occupies an adjacent vertex, they move according to Table \ref{tab:tab3}. Otherwise, they move to the vertex with the smallest index that is not adjacent to the cop's vertex.
		
		\begin{table}[ht!]
			\centering
			\begin{tabular}{|c|c|}
				\hline
				Cop's move & Robber's move \\ \hline
				$i-4$      & $i+1$         \\ \hline
				$i-3$      & $i-1$         \\ \hline
				$i-1$      & $i-3$         \\ \hline
				$i+1$      & $i-4$         \\ \hline
				$i+3$      & $i-4$         \\ \hline
				$i+4$      & $i-4$         \\ \hline
			\end{tabular}
			\caption{Robber's moves in relation to the cop's. The first column corresponds to the vertex the cop just moved to, and the second column is the vertex robber moves to from $i$ as a result of the cop's move.}
			\label{tab:tab3}
		\end{table}
		
		Without loss of generality, the remaining robbers are $R_1$ and $R_2$ and the third robber was captured in round $n-4$. After the third robber was captured, the cop occupies one of $n$, $n-1$ or $n-3$. If she occupies $n$, the remaining robbers move to $n-5$ and $n-7$ respectively, otherwise they move to $n-3$ and $n-7$ in the case of $n-1$ or to $n-1$ and $n-5$ in the case of $n-3$.
		
		We claim that, with the robbers using their described strategy, the cop cannot prevent either from occupying a vertex within [6] until $C^t \in Z$.
		
		If the cop is outside $Z$ and the robber occupies a vertex with a lower index than the cop, we can easily see from Table \ref{tab:tab3} and the robber's strategy of moving to lower indices when possible, that he can move closer to $Z$ if not already there. This is the case when the cop captures the robber at $n$. If the cop then occupies a vertex in $[6]\setminus Z$ and the robber is in $Z$, the robber is able to move to a vertex within $[6]$ at distance 2 from the cop: if the cop occupies 3 or 5, the robber moves to 1 and 2 respectively and if the cop occupies 4 or 6, the robber moves to 6 and 4 respectively. So the robber can remain within $[6]$ until the cop moves to $Z$.
		
		Next we prove that the robbers can return to $[6]$ when the cop captures the robber at $n-3$ or $n-1$. After the robbers move in round $n-4$, it holds that $|C^{n-4}-R_i^{n-4}|=2 \pmod{4}$ for $i\in [2]$. Without loss of generality, we can consider only robber $R_1$. Depending on the vertex that the robber occupies, the cop can move to one of $R_1^{n-4}-3$, $R_1^{n-4}-1$, $R_1^{n-4}+1$, $R_1^{n-4}+3$ or to a vertex at distance at least 2 from the robber. Again considering the robber's strategy, we see that the robber may lower his index regardless of the two vertices he can be at (depending on where the cop captured the first robber) and in particular, he can maintain the difference $|C^t-R_1^t|=2 \pmod{4}$ with his moves by induction on $t\geq n-4$. This remains true until the cop aligns her vertex's index to the robber's modulo $4$, which can only happen once she is in $Z$, since $1$ and $2$ are the unique vertices of the graph with the property $N[1]\cong \mathbb{Z}_4$ (and likewise for $2$) that allows the cop to align her vertex. Since, at the time of capture of the first robber, $d(C^{n-4},Z)>d(R_1^{n-4},[6])$ and $d(C^{n-4},Z)>d(R_2^{n-4},[6])$, the robbers make it back to $[6]$ before the cop makes it to $Z$. Crucially, due to the aforementioned distances, the robbers are within $[6]$ in the round the cop moves to $Z$, otherwise the index of their vertex would have increased before that round.

		The cop requires at least $d_Z=d(n,Z)=d(n-1,Z)=d(n-3,Z)=\ell$ rounds to get to $Z$. Once there, the robbers again avoid capture until they are forced into vertex $n$, which takes another $n-4$ rounds because they were both in vertices in the set $[6]$ when the cop returned to $Z$. In vertex $n$, they split up between $n$ and $n-1$ and again at most one of them is captured in the next round. The last robber repeats the strategy from before by returning to a vertex in $[6]$ and then avoiding capture for $n-4$ more rounds after the cop returns to $Z$.
		With this strategy, the last robber avoids capture until at least round $(n-4)+\ell+(n-4)+\ell+(n-4)=2\ell+3(n-4)$. By Theorem \ref{thm:general_bound} this is also the upper bound for the capture time, so the equality holds.
	\end{proof}
	
	From the proof above, we immediately also obtain the following:
	
	\begin{corollary}
		For $n=4\ell+2$ and $\ell\geq2$, the capture time $\capt(H(n),2)$ is $\ell+2(n-4)$.
	\end{corollary}

	\subsection{2-dismantlable graphs}
	We recall that for 2-dismantlable graphs of order $n$, the capture time is at most $\left\lfloor n/2\right\rfloor$. Using this, we also want to improve the bound for generalized capture time on these graphs. We define the set of vertices $W$ as all vertices in a 2-dismantlable graph, from which the cop can capture the robber in at most $\left\lfloor n/2\right\rfloor$ rounds.
	
	\begin{theorem}
		If $G$ is a 2-dismantlable graph of order $n$ and  
		\begin{equation*}
			d_W=\max\{d(x,W)\mid x\in N[u]\text{, $u$ is a trap}\},
		\end{equation*}
		then $\capt(G,m)\leq m\left\lfloor n/2\right\rfloor+(m-1)d_W$.
	\end{theorem}
	
	\begin{proof}
		The cop starts in $W$. She captures the first robber in at most $\left\lfloor n/2\right\rfloor$ rounds. She then needs at most $d_W$ rounds to return to $W$, from where she can capture the next robber in $\left\lfloor n/2\right\rfloor$ rounds. Repeating this strategy until the end, she requires at most $m\left\lfloor n/2\right\rfloor+(m-1)d_W$ rounds to capture all robbers.
	\end{proof}
	
	For certain families of 2-dismantlable graphs, the upper bound is better still.
	
	\begin{theorem}
		\label{thm:upper_paths}
		\begin{enumerate}
			\item For $n\geq1$ and $m\geq 2$, $\capt(P_n,m)=n-1$.
			\item For $n\geq2$ and $m\leq n-1$, $\capt(K_n,m)=m$.
			\item For $n\geq4$ and $m\leq n-1$, $\capt(W_n,m)\leq2m-1$.
		\end{enumerate}
	\end{theorem}
	
	\begin{proof}
		For paths, we label the vertices from 1 to $n$ so that vertices 1 and $n$ are leaves and each vertex $i$ is adjacent to $i-1$ and $i+1$ for $2\leq i\leq n-1$. The cop wins in round $n-1$ if she starts at vertex 1 and moves from $i$ to $i+1$ in each round. The robbers force the length of the game to be $n-1$ by having a single robber start at vertex 1 and the others at vertex $n$, then they skip all rounds until the end.
		
		For complete graphs, the cop starts at any vertex. In every round, she is adjacent to at least one robber, whom she captures next. She requires at most $m$ rounds to capture all robbers. Conversely, the robbers pick $m$ distinct vertices that are not the cop's vertex and then skip every round, ensuring it takes $m$ rounds to capture them all.
		
		For wheels, the cop starts at the universal vertex. She is adjacent to at least one robber, whom she captures in round 1. In round 2, she returns to the universal vertex, from where she captures the next robber in round 3. Repeating this, she captures all robbers in round $2m-1$ at the latest.
	\end{proof}
	
	\begin{remark}
		The capture time of complete graphs and wheels is more difficult to determine when $m\geq n$. In Propositions \ref{prop:complete} and \ref{prop:wheels} we provide a limit as the number of robbers continues to increase for both graph classes.
	\end{remark}
	
	We know due to \cite{MEHRABIAN2011102} that the capture time of a tree $G$ is at most $\left\lceil\diam(G)/2\right\rceil$. We have the following result for more robbers.
	
	\begin{theorem}
		\label{thm:trees}
		Let $G$ be a tree of order $n$ with $\ell$ leaves and let $m\leq \ell$. Then $\capt(G,m)\leq \left\lceil\diam(G)/2\right\rceil+(m-1)\diam(G)$.
	\end{theorem}
	
	\begin{proof}
		The cop starts in the center of the tree. From there, she captures the first robber in $\left\lceil\diam(G)/2\right\rceil$ rounds. She captures each remaining robber in at most $\diam(G)$ rounds, giving us the bound.
	\end{proof}
	
	This bound is tight. To conclude the section, we provide an infinite family of trees for which the bound is reached for all $m$.
	
	\begin{proposition}
		\label{thm:trees_tight}
		If $m\geq 1$ and $G$ is the subdivided star $T_{m+1}$, then $\capt(G,m) = 4m-2$.
	\end{proposition}
	
	\begin{proof}
		As $\diam(G) = 4$, by Theorem \ref{thm:trees}, we only need to prove that the robbers can survive at least $4m-2$ rounds.
		They start by picking distinct vertices in $G$ which are all leaves and they are all as far from the cop's starting vertex as possible. This ensures that the distance between the cop and any robber is at least 2 in round 0. Then they skip every round, requiring the cop to reach the upper bound.
	\end{proof}
	
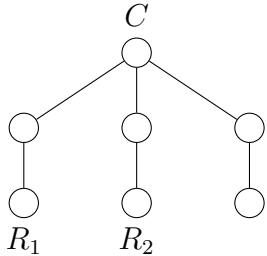
\begin{figure}[ht!]
		\centering
		\begin{tikzpicture}[node/.style ={draw, circle}]
			\node[node] (1) at (0,0) {};
			\node[node] (2) at (-1.5,-1) {};
			\node[node] (3) at (0,-1) {};
			\node[node] (4) at (1.5,-1) {};
			\node[node] (5) at (-1.5,-2) {};
			\node[node] (6) at (0,-2) {};
			\node[node] (7) at (1.5,-2) {};
			\draw (1) -- (2);
			\draw (1) -- (3);
			\draw (1) -- (4);
			\draw (2) -- (5);
			\draw (3) -- (6);
			\draw (4) -- (7);
			\node (A) at (0,0.5) {$C$};
			\node (B) at (-1.5,-2.5) {$R_1$};
			\node (C) at (0,-2.5) {$R_2$};
		\end{tikzpicture}
		\caption{An example of the tree attaining the bound in Theorem \ref{thm:trees}.}
		\label{fig:tree_max}
	\end{figure}

	For example, see the subdivided star $T_3$ ($m=2$) and the starting positions of the robbers in Figure \ref{fig:tree_max}.
	
	\section{Increasing the number of robbers}
	We notice that for paths, the game has a finite length regardless of the number of robbers we add. However, this does not hold for all cop-win graphs. In this section, we expand this idea further. As we continue to add robbers, we want to find out on which graphs the cops can capture an infinite number of robbers in finite time. More precisely, we wonder for which graphs the limit $\lim_{m\rightarrow \infty} \capt_k(G,m)$ exists.
	
	\begin{definition}
		If $\lim_{m\rightarrow \infty} \capt_k(G,m)$ exists, we say that $G$ is a {\em strong $k$-cop-win graph}. If $k=1$, we simply say that $G$ is a {\em strong cop-win graph}.
	\end{definition}
	
	Observe that every graph $G$ on $n$ vertices is strong $n$-cop-win by definition. Thus a more relevant question is to find the smallest $k$ such that $G$ is strong $k$-cop-win. 
	
	\begin{remark}
		\begin{itemize}
			\item If $G$ is strong $k$-cop-win, then it is also strong $\ell$-cop-win for every $\ell \geq k$.
			\item If $G$ is not $k$-cop-win, then it is also not strong $k$-cop-win.
		\end{itemize}
	\end{remark}
	
	\begin{proposition}
		\label{prop:limits-examples}
		If $n\geq 2$, the following graphs are strong cop-win with limits:
		\begin{enumerate}
			\item $\lim_{m\rightarrow \infty} \capt(P_n,m)=n-1$ for paths,
			\item $\lim_{m\rightarrow \infty} \capt(K_{1,n},m) = 2n-1$ for stars and
			\item $\lim_{m\rightarrow \infty} \capt(G,m) = n + \ell - 2$ for caterpillars with $\ell$ leaves and $n$ vertices.
			
		\end{enumerate}
	\end{proposition}
	
	\begin{proof}
		The result for paths follows from the proof of Theorem \ref{thm:upper_paths}.
		
		For stars, we find a strategy for the cop that captures an arbitrary number of robbers in the same number of rounds. For a star $K_{1,n}$, we label the universal vertex as $u$ and we number the other vertices from 1 to $n$. The cop starts in $u$. She repeats the following moves. If she occupies $u$, she moves to $i$, $1\leq i \leq n$, where she has not occupied $i$ earlier in the game. If she occupies a leaf, she moves to $u$. We need to prove that once she has visited all leaves, all the robbers have been captured. This is equivalent to proving that a robber who starts in some vertex $i$ cannot move to any other vertex $j\neq i$ without being captured in the next round. Since he cannot move to a vertex that the cop has already checked, he is captured once the cop checks $i$. Our claim is obvious, since a robber can only move from $i$ to $u$, which the cop will occupy in the next round at the latest. For all $m\geq n$, the robbers can occupy all the leaves of the graph at the beginning of the game and then skip every turn, forcing the game to last until round $2n-1$ and proving the limit.
		
		Lastly, we need to prove the claim for caterpillars. The cop follows the central path from one end-vertex to another in a way in which she checks all leaves adjacent to a vertex of the central path in the same manner as in the star graph. The same proof as before shows that a robber in a leaf, adjacent to the vertex on the central path that the cop has already visited, cannot reach any vertex that the cop has checked before without being captured in the next round and no robber on the central path can reach a vertex on the central path which was already checked by the cop without being captured in the next round. Therefore, the cop captures all robbers once she checks all vertices of the caterpillar at least once. This takes her at most $n-\ell-1$ moves to check the central path and $2\ell -1$ moves to check the leaves, excluding the last. A sufficient number of robbers can occupy all leaves at the start, forcing $n+\ell -2$ to be the length of the game.
	\end{proof}
	
	It is well-known that one cop cannot win the game on a cycle $C_n$ for $n \geq 4$. However, we have the following result.
	
	\begin{proposition}
		Cycles are strong 2-cop-win graphs with the limit $\lim_{m\rightarrow \infty} \capt_2(C_n,m)=\lfloor \frac{n-1}{2} \rfloor$.
	\end{proposition}
	
	\begin{proof}
		The cops start the game on adjacent vertices of the cycle. In subsequent rounds, they move around the cycle in different directions, each taking the shortest path to the opposite side of the cycle. The robbers cannot reach the vertices the cops have already been to, so the game ends at the latest in round $\lfloor \frac{n-1}{2} \rfloor$, since the cops are either at two adjacent vertices when $n$ is even or at the same vertex when $n$ is odd and have checked every vertex. The robbers simply have to stay at the vertex farthest from one of the cops in the above case. If the cops start on non-adjacent vertices, there are two internally disjoint paths between their starting vertices. If there are at least two robbers, one of the robbers starts on one of the paths, while the others start on the other, and all start so that they are as far from the cops as possible. The cops have to capture either the lone robber or the remaining group first, and then have to spend as many moves as in the original case for the remaining robbers. If cops try to capture all the robbers at the same time by moving in the same direction, the robbers can evade capture forever.
	\end{proof}
	
	The result for caterpillars proves stronger. But first we prove a general result about the behavior of strong $k$-cop-win graphs under retracts.
	
	\begin{lemma}
		\label{lem:retract}
		If $H$ is a 1-point-retract of $G$ and $m, k\geq 1$, then $\capt_k(H,m) \leq \capt_k(G,m)$.
	\end{lemma}
	
	\begin{proof}
		Let $H$ be a 1-point-retract of $G$ with the 1-point-retraction $f$ and let $\capt_k(G,m) = \ell$. We describe a strategy of $k$ cops that can catch $m$ robbers on $H$ in at most $\ell$ rounds.
		
		While the game with $k$ cops and $m$ robbers is played on $H$, we imagine a parallel game on $G$ (also with $k$ cops and $m$ robbers). Robbers are playing optimally on $H$, and we simply copy their moves to $G$ (which is legal as $f$ is an identity on $H \subseteq G$). The cops are playing optimally on $G$ and can catch all robbers in at most $\ell$ rounds. When cop $c$ moves from $x$ to $y$ in $G$, we let $c$ move from $f(x)$ to $f(y)$ in $H$ (which is legal as $f$ is a homomorphism). If cop $c$ captures a robber $r$ in $G$, so $c=r$, then by our assumption that $r \in V(H)$, we know that $f(c)=c$, thus cop $c$ also captures robber $r$ in the game on $H$. This shows that by the time the last robber is captured in the game on $G$, the last robber is also already captured in the game on $H$.
	\end{proof}
	
	\begin{proposition}
		\label{thm:caterpillars-strong}
		A tree $T$ is strong cop-win if and only if $T$ is a caterpillar.
	\end{proposition}
	
	\begin{proof}
		Caterpillars are strong cop-win by Proposition \ref{prop:limits-examples}. Each tree that is not a caterpillar includes the subdivided star $T_3$ as an induced subgraph. Denote the vertex of degree 3 by $u$. We begin by proving that this tree is not strong cop-win.
		
		We need to prove that $\lim_{m\rightarrow \infty} \capt_k(T_3,m)$ diverges. Without loss of generality, let $m=2^\ell$ for some $\ell\geq2$. We prove using induction on $\ell$ that, if the robbers start in the leaves that are as far as possible from the cop's initial vertex, then $\capt_k(T_3,m)\geq4\ell+\ecc(C^0)$, where $C^0$ is the cop's starting vertex. The base case $\ell=1$ is just the proof of Theorem \ref{thm:trees_tight}. Assume that the theorem holds up to $\ell$ and consider $\ell+1$. The cop starts in some arbitrary vertex in $T_3$. The robbers split into two groups of $2^\ell$ robbers. Each group starts in a distinct leaf, picking the two leaves furthest from the cop. This way, the cop requires at least $\ecc(C^0)$ rounds to capture one group. The robbers skip every round until the cop is adjacent to one group. That group skips, while the other group splits into two -- $2^{\ell-1}$ robbers move to the adjacent vertex, while the other $2^{\ell-1}$ skip. We may assume that the cop will move to capture the adjacent group, otherwise the groups switch roles. The group of robbers that previously moved, now moves to $u$, while the other remaining robbers skip. The cop cannot capture any robbers in the next round, so the group of robbers occupying $u$ moves to the vertex adjacent to the leaf that neither the cop nor any robber occupied until then. In the next round, this group moves to the leaf, while the cop is somewhere between the leaf in which she captured the first group and $u$. In any case, $\ecc(C)\geq 2$ and the robbers occupy the two leaves furthest from the cop. Using the induction hypothesis, the cop requires at least $4\ell+2$ rounds to win, meaning she needs at least $\ecc(C^0)+2+4\ell+2=\ecc(C^0)+4(\ell+1)$ rounds in total.
		As $f(\ell)=4\ell$ is a monotonically increasing function that bounds the capture time from below, the limit diverges.
		
		If $T$ is a tree that is not a caterpillar, $T_3$ is its retract. By Lemma \ref{lem:retract}, the capture time of $T$ is greater than that of $T_3$ for any number of robbers, so it is not a strong cop-win graph.
	\end{proof}
	
	\begin{theorem}
		Graphs $H(n)$, where $n=4\ell+2$ for $\ell\geq 2$, are not strong cop-win.
	\end{theorem}
	
	\begin{proof}
		We denote $t=\ell+n-4$ and let $m=3^s$ for $s\geq1$. We will prove by induction that, if all robbers start in $[6]$ and none are adjacent to the cop, $\capt(G,m)\geq st$. The case $s=1$ is simply the proof of Theorem \ref{thm:tight_hn}. Assume the theorem hold until $s$ and we prove for $s+1$. The robbers follow a modified strategy from the proof of Theorem \ref {thm:tight_hn}. They all start in the same vertex and follow their strategy until they are forced into the trap. There, they split up into equal thirds, meaning at least $\left\lfloor2/3\cdot 3^{s+1}\right\rfloor=2/3\cdot 3^{s+1}=2\cdot 3^s>3^s$ robbers remain uncaptured after $n-4$ rounds. All of these robbers follow the strategy laid out in  the proof of Theorem \ref {thm:tight_hn}, so that the cop is forced to return to $Z$ before capturing another robber. This takes another $\ell$ rounds. Since the robbers followed their strategy, at least $3^s$ return to $[6]$ and they are not adjacent to the cop, so we may use the induction hypothesis. With this, the game lasts at least $n-4+\ell+st=t+st=(s+1)t$ rounds. The function $f(s)=s$ is a monotonically increasing, non-constant function that bounds the capture time from below, proving the theorem.
	\end{proof}
	
	\begin{proposition}
		\label{prop:complete}
		Complete graphs are not strong $k$-cop-win for $1\leq k\leq \left\lceil n/2\right\rceil-1$ and are strong $\left\lceil n/2\right\rceil$-cop-win. Moreover, $\lim_{m \to \infty} \capt_{\lceil n/2 \rceil} (K_n, m) = 1$.
	\end{proposition}
	
	\begin{proof}
		Let $1\leq k\leq \left\lceil n/2\right\rceil-1$. We prove that $\capt_k(K_n,m)\geq \ell$ by induction on $\ell$.
		The case $\ell=1$ is just the proof of Theorem \ref{thm:upper_paths}. We assume the theorem holds for $\ell$ and consider $\ell + 1$. 
		It is enough to prove the claim for $k=\left\lceil n/2\right\rceil-1 = \left\lceil (n-2)/2\right\rceil$ as it then follows immediately for any smaller $k$. Let $m=(n- \left\lceil (n-2)/2\right\rceil)^{\ell+1}=\left\lceil (n+1)/2\right\rceil^{\ell+1}$. After the cops choose any starting vertices in round 0, the robbers choose their starting vertices as all unoccupied vertices so that there are $\left\lceil (n+1)/2\right\rceil^\ell$ robbers at each vertex they occupy. In round 1, cops can move to capture at most $\left\lceil (n-2)/2\right\rceil \left\lceil (n+1)/2\right\rceil^\ell$ robbers, and since $n \geq 2\left\lceil (n-2)/2\right\rceil+1$, at least $(n-2\left\lceil (n-2)/2\right\rceil) \left\lceil (n+1)/2\right\rceil^\ell \geq \left\lceil (n+1)/2\right\rceil^\ell$ robbers remain in the game. These robbers can move anywhere in round 1 so we can use the induction hypothesis. This implies that $\capt_k(K_n, m) \geq \ell$ and as $f(\ell)=\ell$ is a monotonically increasing function, the limit $\lim_{m\to \infty} \capt_k(K_n, m)$ does not exist.
		
		For $k=\left\lceil n/2\right\rceil$, the cops start on distinct vertices. The robbers can start anywhere on the remaining $\left\lfloor n/2\right\rfloor$ vertices. The cops then move to capture all of them in round 1, proving that $\lim_{m\rightarrow \infty} \capt_k(K_n,m) = 1$.
	\end{proof}

	\begin{proposition}
		\label{prop:wheels}
		Wheels are not strong cop-win but are strong 2-cop-win. Moreover, 
		\begin{equation*}
			\lim_{m \to \infty} \capt_2(W_n, m) = n-3.
		\end{equation*}
	\end{proposition}
	
	\begin{proof}
		Let $u$ be the universal vertex of $W_n$ and let $1, \ldots, n-1$ be the vertices on the cycle. We first prove that $\capt_2(W_n, m) \leq n-3$ for every $m \geq 1$. One cop starts the game on $u$ and then alternates between $u$ and $n-1$ for the rest of the game (thus catching any robber entering $u$ or $n-1$ during the game). The second cop starts the game on $1$ and moves to $2, 3, \ldots, n-2$ consecutively, capturing any robber positioned on these vertices. Thus all robbers are caught by round $n-3$ at the latest. 
		
		Now we need to prove that the robbers can evade capture until then. Put the vertices into two classes. If a player had already visited a vertex and the last player to visit was a cop, the vertex is {\em checked}, otherwise it is {\em unchecked}. Consider the game with infinite robbers. At the start of the game, the cops start at any two vertices. Meanwhile, the robbers start with infinite robbers at each vertex, so $n-2$ vertices are unchecked. Clearly, the cops must reduce the number of unchecked vertices to 0 in order to win the game. We notice that with the cops' strategy above, the number of unchecked vertices strictly decreases after each round and reaches 0 in round $n-3$. Now consider an arbitrary strategy of two cops and let the robbers play according to the following strategy: in each round, each infinite group of robbers at a vertex considers their neighboring vertices and occupies each unoccupied neighbor so that the original vertex and its neighbors are each occupied by infinite robbers. We claim that the number of checked vertices increases by at most one after each round (i.e.\ after both sets of players move), except maybe in the last round. If the cops move so that they only add one new checked vertex, there is nothing to prove. Now suppose that the cops move such that they add two new checked vertices. This means that both cops move to a new vertex, leaving two previously checked vertices unoccupied. While the game is not over yet and the robbers are using the described strategy, at least one of these two previously checked vertices is adjacent to a vertex occupied by the robbers as the graph is 2-connected. By the robbers' strategy, they now change the status of this vertex from checked to unchecked, thus the number of checked vertices can increase by at most one after each round. After round $n-4$, at most two vertices are checked, thus all $n$ vertices are checked in round $n-3$ at the earliest (after the cops move, capturing the last robbers and so the robbers cannot move in this round). Clearly, the robbers' strategy still works with finitely many robbers, as long as there are enough.
		
		It remains to prove that wheels are not strong cop-win. Let $m = (n-1) 2^\ell$, $\ell \geq 0$. We prove that $\capt(W_n, (n-1) 2^\ell) \geq \ell$. This trivially holds for $\ell = 0$. Now let $\ell \geq 1$. In round 0, robbers position themselves on all vertices not occupied by the cop, $2^\ell$ on each. We can assume that cop moves in round 1, and thus captures $2^\ell$ robbers. The now unoccupied vertex has at least two neighbors occupied by robbers, and one of them moves $2^{\ell-1}$ robbers to the empty vertex. Now every vertex not occupied by the cop contains at least $2^{\ell-1}$ robbers, thus by the induction hypothesis, $\capt(W_n, (n-1) 2^\ell) \geq 1 + \capt(W_n, (n-1) 2^{\ell-1}) \geq 1 + (\ell-1) = \ell$. As $f(\ell) = \ell$ is a monotonically increasing function, the limit $\lim_{m \to \infty} \capt(W_n, m)$ does not exist.
	\end{proof}
	
	\begin{proposition}
		\label{prop:complete-bipartite}
		If $n \geq m \geq 1$, then $K_{m,n}$ are strong $m$-cop-win graphs and are not strong $k$-cop-win graphs for $1 \leq k \leq m-1$. Furthermore, 
		\begin{equation*}
			\lim_{\ell \to \infty} \capt_m(K_{m,n}, \ell) = \begin{cases}
				2\left\lceil\frac{n}{m}\right\rceil-2; & n > m\\
				1; & n = m.
			\end{cases}
		\end{equation*}
	\end{proposition}
	
	\begin{proof}
		Let $(M, N)$ be the bipartition of $K_{m,n}$ where $|M| = m$ and $|N| = n$. In the game with $m$ cops, their strategy is to start the game by occupying all vertices of $M$, and then visit different subsets of $m$ vertices in $N$ each odd round, returning to $M$ in the even rounds. Thus robbers cannot move between vertices without being captured and thus $\capt_m(K_{m,n}, t) \leq \left \lceil \frac{n}{m} \right \rceil$. As the sequence $(\capt_m(K_{m,n}, t))_{t \geq 1}$ is an increasing integer sequence, this already implies that $\lim_{t \to \infty} \capt_m(K_{m,n}, t)$ exists and that $K_{m,n}$ is strong $m$-cop-win.
		
		Suppose now that $k \leq m-1$ cops are playing the game. We prove that $\capt_k(K_{m,n}, (m+n) n^\ell) \geq \ell$, which implies that the limit $\lim_{t \to \infty} \capt_k(K_{m,n}, t)$ diverges. The bound is trivially true for $\ell = 1$. Let $\ell \geq 2$. The robbers' strategy is to place $n^\ell$ robbers on each vertex of the graph. After round 1, as $k \leq m-1$, at least one vertex in $M$ and at least one in $N$ are still occupied by the robbers. From the one in $M$, $n^{\ell-1}$ robbers are moved to each vertex in $N$, and from the one in $N$, $n^{\ell-1}$ are placed on each vertex in $M$. Now there are at least $n^{\ell-1}$ robbers on each vertex in in the graph (that is not occupied by a cop). Using the induction hypothesis, we obtain $\capt_k(K_{m,n}, (m+n) n^\ell) \geq 1 + \capt_k(K_{m,n}, (m+n) n^{\ell-1}) \geq \ell$. 
		
		For capture time in the case $m=n$, it is enough for cops to start in any $m$ distinct vertices and move to capture all robbers in round 1. Conversely, the robbers can start distributed on all unoccupied vertices in round 0 to force capture time at least 1.
		
		Let us now consider the capture time in the case $n > m$. Let the game be played by $m$ cops and $\ell$ robbers. The cops can win in at most $2\lceil n/m\rceil-2$ rounds by starting in $m$ distinct vertices of $N$, then moving to all $m$ vertices of $M$ in odd rounds and $m$ distinct, previously unvisited vertices of $N$, in even rounds. As cops occupy all vertices of $M$ in every second round, the robbers that move between vertices are captured immediately. Thus the above strategy of the cops ensures that the game ends at the latest in round $2\lceil n/m\rceil-2$.
		
		Still considering the case $n > m$ it remains to prove that the capture time is not less than $2\lceil n/m\rceil-2$ (if there are enough robbers playing the game). Let the game be played by $m$ cops and $\ell \geq n+m$ robbers. Robbers' strategy is to start the game by placing one robber on each vertex and then never moving during the game. We use checked and unchecked vertices as defined in the proof of Proposition \ref{prop:wheels}. If cops want to increase the number of checked vertices, they need to alternate between $M$ and $N$. After round 2, at most $3m$ vertices are checked, where at most $m$ of these vertices are from $M$ and at most $2m$ from $N$. From now on we only consider vertices from $N$. After every next two rounds, at most $m$ new vertices become checked in $N$. Thus (just to achieve that all vertices from $N$ are checked) the number of rounds is at least $2 + 2 \lceil \frac{n-2m}{m} \rceil = 2 \lceil \frac{n}{m} \rceil - 2$.
	\end{proof}
	
	We notice that so far, the 0-visibility cop number and the smallest number of cops required for the graph to be strong $k$-cop-win match for all considered families of graphs. The next result proves that they are always equal.
	We call the maximum index of the round in which the robber is captured by the cops in the 0-visibility game on graph $G$ the \emph{0-visibility capture time of $G$}.
	
	\begin{theorem}
		\label{thm:characterization}
		It holds that $c_0(G) \leq k$ if and only if $\lim_{m\rightarrow \infty} \capt_k(G,m) < \infty$. Furthermore, the 0-visibility capture time and the limit are equal if both exist.
	\end{theorem}
	
	\begin{proof}
		As $k$ cops win the 0-visibility game played on $G$ no matter how robber plays, let $\ell$ be the 0-visibility capture time of $G$. Let the $k$ cops play according to their winning strategy in the 0-visibility game on $G$ against $m$ cops (cops can imagine that robbers are invisible). If after $\ell$ moves, some robber $r$ has not been caught yet, this contradicts the fact that cops were using their winning strategy in the 0-visibility game on $G$ (a robber in that game could have taken $r$'s trajectory). Thus for all $m \geq 1$, $\capt_k(G,m) \leq \ell$ and $\lim_{m\rightarrow \infty} \capt_k(G,m) < \infty$.
		
		Conversely, assume that $\lim_{m\rightarrow \infty} \capt_k(G,m) = \ell$. Let $\mathcal{M}$ be the set of all possible sequences $(v_i)_{0\leq i\leq \ell}$ of vertices in $G$, for which it holds that $v_{i+1}$ and $v_i$ are adjacent (note that since the graph is reflexive, each vertex is adjacent to itself) and let $m=|\mathcal{M}|$. Since the limit exists and the sequence $(\capt_k(G, m))_{m \geq 1}$ is increasing, $\capt_k(G,m)\leq \ell$ and we claim that equality holds. Every game played by $m$ robbers is finite, including the game in which each robber uses exactly one strategy from $\mathcal{M}$. If equality doesn't hold, the $m$ robbers are captured before round $\ell$. Since the limit equals $\ell$, there is $m'>m$ for which $\capt_k(G,m')= \ell$. However, the robbers in the first case are already using all distinct strategies, so $m'-m$ robbers are using the same strategy as some other robber. Since the remaining robbers are captured in fewer than $\ell$ moves, so are the $m'-m$ robbers playing the exact same strategy, as they are captured in the same move as their copies. Thus all robbers are captured in fewer than $\ell$ moves, contradicting the assumption.
		
		We call the $m$-tuple $(R_1^i, R_2^i,\ldots, R_m^i)$ the {\em configuration} of robbers after round $i$. We claim that the strategy used by the cops against the robbers in the game with $m$ robbers playing all possible distinct strategies is also the winning strategy in the 0-visibility game. Since the robbers are using all possible distinct strategies of length $\ell$ at the same time, the configuration of robbers is unique up to their indexing in each round of the game, meaning the cops know all of configurations in advance. If there existed any other configuration of robbers after a particular round, at least two robbers would be playing the same strategy because we assumed that we have enough robbers so that they are able to play each distinct strategy on the graph. But if any two robbers are playing the same strategy, the cops may continue to use the same strategy they would otherwise, for they will capture both of those robbers when they would otherwise capture one. The cops do have a winning strategy, meaning that this strategy is independent of the configuration of the robbers in each round. Now, the cops need not see the robbers, but simply pretend to play against these $m$ robbers in the 0-visibility game, where one of these $m$ robbers is real, while the others are not. Using their strategy, the cops capture this robber when they would capture him in the original game. Since the cops win against all distinct strategies of the robber, it does not matter how the robber plays in the 0-visibility game. Thus $c_0(G) \leq k$ and as the robber can play according to the strategy from $\mathcal{M}$ that is caught only in round $\ell$, the 0-visibility capture time is exactly $\ell$.
		
		The robber in the 0-visibility game may play the strategy of the robber that is captured last in the regular game to show that the 0-visibility capture time is at most $\ell$ if the limit exists. 
		We already know from the first part of the proof that if the 0-visibility capture time is $\ell$, then $\lim_{m\rightarrow \infty} \capt_k(G,m) \leq \ell$. Since the existence of one implies the existence of the other, they are equal.
	\end{proof}
	
\begin{figure}[ht!]
		\centering
		\begin{tikzpicture}[node/.style ={draw, circle}]
			\node[node] (1) at (0,0) {};
			\node[node] (2) at (-1.5,-1) {};
			\node[node] (3) at (-1.6,0.7) {};
			\node[node] (4) at (0,1.7) {};
			\node[node] (5) at (1.6,0.7) {};
			\node[node] (6) at (1.5,-1) {};
			\draw (1) -- (2);
			\draw (1) -- (3);
			\draw (1) -- (4);
			\draw (1) -- (5);
			\draw (1) -- (6);
			\node (A) at (0,-1) {$\cdots$};
		\end{tikzpicture}
		\caption{The infinite star.}
		\label{fig:star_inf}
	\end{figure}
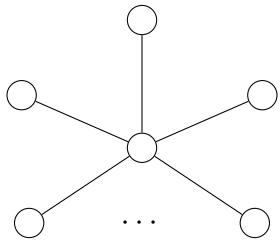

	The theorem applies only to finite graphs. Consider the following example.
	\begin{example}
		Let $G$ be the infinite star $K_{1,\infty}$ (see Figure \ref{fig:star_inf}). Clearly, $\lim_{m\rightarrow \infty} \capt_k(G,m)$ diverges. In the 0-visibility game, however, the cop can choose an ordering of leaves and win by starting at the universal vertex $u$ and then moving to the next vertex in the ordering and the universal vertex interchangeably. The robber chooses his starting vertex and gets captured either when the cop checks his vertex or, if he moves to the universal vertex, in the next round.
	\end{example}
	
	\begin{remark}
		The proof of Theorem \ref{thm:characterization} reveals even more about strong cop-win graphs. That is, the cops in strong cop-win graphs can always play the same strategy against any number of robbers, and that strategy is independent of the robbers' moves. We may call this strategy \emph{a strong winning strategy for the cops}. This strategy is not the fastest in general -- the proof of Theorem \ref{thm:upper_paths} shows that the cop wins against any number of robbers on a path by checking the path from one leaf to the other. However, she can clearly win faster against one robber as paths are 2-dismantlable. The strategy is the fastest from a certain number of robbers onward, with the maximum given by $m$ from the proof of Theorem \ref{thm:characterization}.
	\end{remark}
	
	\section{Further directions}
	
	In this paper we define a variant of the game of cops and robbers with multiple robbers and study its capture time. However, there are several directions left that we have not explored. For example, it would be interesting to obtain sharp general upper bounds for the capture time of of $k$-cop-win graphs for $k \geq 2$, and to consider possible improvements of the bound from Theorem \ref{thm:general_bound} for cop-win graphs.
	
	Theorem \ref{thm:characterization} yields a surprising connection between the generalized capture time and the 0-visibility cops and robber. It would be worth exploring $\lim_{m\rightarrow \infty} \capt_k(G,m)$ further as it could give useful insight into the open problems in the area of limited visibility cops and robber.
	
	We also wonder if there is an interesting relation between the capture time in the game with multiple robbers as defined in \cite{BONATO20095588} and the capture time defined in this paper.
	
	\section*{Acknowledgments}
	
	V.I.C.\ acknowledges the financial support from the Slovenian Research and Innovation Agency (Z1-50003, P1-0297, N1-0218, N1-0285, N1-0355) and the European Union (ERC, KARST, 101071836).

\end{document}